\documentclass[10pt]{article}
\usepackage{amsmath,amssymb,amsthm,amscd}
\numberwithin{equation}{section}

\newtheorem{prop}{Proposition}[section]
\newtheorem{theo}[prop]{Theorem}
\newtheorem{lemm}[prop]{Lemma}

\newtheorem{rema}[prop]{Remark}

\newtheorem{conj}[prop]{Conjecture}
\newtheorem{ack}[prop]{Acknowledgment}
\def\begeq{\begin{equation}}
\def\endeq{\end{equation}}

\def\and{\quad{\rm and}\quad}

\def\<{\langle}
\def\>{\rangle}

\def\lab{\label}

\begin{document}

\title{a Modified K\"ahler-Ricci Flow}

\author{Zhou Zhang\\
Department of Mathematics\\
University of Michigan, at Ann Arbor}
\date{} 
\maketitle

Abstract: In this note, we study a K\"ahler-Ricci 
flow modified from the classic version. In the 
non-degenerate case, strong convergence at 
infinity is achieved. We also have partial results 
for some interesting degenerate cases.     

\section{Set-up and Motivation}

K\"ahler-Ricci flow, which is nothing but Ricci flow with initial metric being K\"ahler, 
enjoys the same debut as Ricci flow in R. Hamilton's original paper \cite{ham}. 
H. D. Cao's paper, \cite{cao}, can be seen as the first one devoted to the study of 
K\"ahler-Ricci flow and the alternative proof of Calabi's conjecture presented there 
has been bringing great interest to this object.  

Though it is essentially Ricci flow, the cohomology meaning coming with K\"ahler 
condition makes it possible to transform it to a equivalent scalar flow \footnote{This  
statement makes use of the uniqueness and short time existence results of Ricci flow.},  
which is much simpler-looking and more flexible to study. The discussion in this note 
would hopefully give a flavor of the flexibility.\\ 

Let $\omega_0$ be any K\"ahler metric over a closed manifold $X$ (with complex 
dimension greater or equal to $2$), and $\omega_\infty$ is any smooth real closed 
$(1,1)$-form. Set $\omega_t=\omega_\infty+e^{-t}(\omega_0-\omega_\infty)$ and 
consider the following flow over the level of metric potential for space-time: 
\begin{equation}
\lab{eq:flow}\frac{\partial u}{\partial t}={\rm log}\frac{(\omega_t+\sqrt{-1}\partial\bar
{\partial}u)^n}{\Omega}, ~~~~u(0,\cdot)=0,
\end{equation}
where $\Omega$ is a smooth volume form over $X$. 

Let $\tilde{\omega}_t=\omega_t+\sqrt{-1}\partial\bar{\partial}u$ and the corresponding 
flow on the level of metric is a little bit artificially looking as follows: 
\begin{equation}
\lab{eq:MKR} \frac{\partial\tilde{\omega}_t}{\partial t}=-{\rm Ric}(\tilde{\omega}_t)+{\rm 
Ric}(\Omega)-e^{-t}(\omega_0-\omega_\infty), ~~\tilde{\omega}_0=\omega_0,
\end{equation}
where the meaning of the form, ${\rm Ric }(\Omega)$, as in \cite{t-znote}, is a natural 
generalization from the Ricci form for a K\"ahler metric, i.e., using the volume form 
$\Omega$ instead of the volume form for some K\"ahler metric in the expression of 
Ricci form from classic computation in K\"ahler geometry.    

\begin{rema}

The equation (\ref{eq:MKR}) doesn't look so natural at the first sight when $\omega_0
\neq\omega_\infty$, but it's essentially still a K\"ahler-Ricci flow, and the extra term in 
comparison to the flow studied in \cite{cao}, which is exponentially decaying, should 
not bring too much difference in spirit. 

Our motivation to study this flow is to solve the following complex Monge-Amp\`ere 
equation 
$$(\omega_\infty+\sqrt{-1}\partial\bar{\partial}u_\infty)^n=\Omega,$$    
using flow techniques. This has been done in the case of $[\omega_\infty]$ being 
K\"ahler in \cite{cao}, which provides another proof of Calabi's conjecture. 

One can also solve it for some degenerate $[\omega_\infty]$ (semi-ample and big) 
by method of continuity using other (more direct) perturbation, which seems to be 
less delicate than K\"ahler-Ricci flow as described in \cite{thesis} and \cite{tosatti}.   

The point is to allow the change of cohomology class along the flow, which is important 
in the consideration of $[\omega_\infty]$ being degenerate as a K\"ahler class. The 
modification of original K\"ahler-Ricci flow by such a term as above is inevitable from 
simple cohomology consideration.   

\end{rema}

Our results can be summarized in the following theorem.

\begin{theo}

The modified K\"ahler-Ricci flow $(1.1)$ (or $(1.2)$ equivalently) exists as long as the 
cohomology class, $[\omega_t]$ remains K\"ahler. 

1) When $[\omega_\infty]$ is K\"ahler, the flow converges exponentially smoothly to the 
unique solution of the corresponding Monge-Amp\`ere equation; 

2) When $[\omega_\infty]$ is semi-ample and big, we have degenerate estimates on the 
metrics along the flow out of the stable base locus set of $[\omega_\infty]$ uniform for all 
time and the volume form, $\tilde\omega_t$, is bounded from above and away from 0 
along the way.  

3) When $[\omega_\infty]$ is "only" big, i.e., the flow exists up to a finite time $T$, and 
$[\omega_T]$ is semi-ample, we have local smooth convergence of the flow out of the 
stable base locus set of $[\omega_T]$.   

\end{theo}

The rest part of this note will be devoted to the proof of this theorem. 

\begin{ack}

I would like to thank my thesis advisor, Professor Gang Tian, for introducing this interesting 
field and encouragement along the way. This research was partially done during the stay 
at MSRI (Mathematical Sciences Research Institute) as a Posdoctor Research Fellow on 
academic leave. I would like to thank the institute and Department of Mathematics of 
University of Michigan, at Ann Arbor, for their kindness and effort to make this opportunity 
possible for me. Also the hospitality of the institute can not be appreciated enough. 

Several people, Yanir Rubinstein and others, took the trouble to read an earlier version of 
this note. I would like to thank them for all their feedbacks.  

\end{ack}

\section{General Facts and Basic Computations}

The equation $(1.1)$ is clearly still parabolic, and so short time existence and uniqueness 
is not a problem. It's also easy to see that the smooth solution exists as long as $[\omega_t]$ 
remains K\"ahler as already being described in \cite{thesis}. Simply speaking, when arguing 
locally in time for this range, $\omega_t$ can be made uniform as metric which makes life 
very easy to follow Cao and Yau's argument as in \cite{cao} and \cite{yau}. So the existence 
part of Theorem 1.2 is justified. \\ 

Convergence, or estimate uniform for time, is our main concern now. For all the expressions 
below, $C$ would be a positive constant (fixed for each place). Let's list some basic 
computation from $(1.1)$ in the following. 
$$\frac{\partial}{\partial t}(\frac{\partial u}{\partial t})=\<\tilde{\omega}_t, \frac{\partial\tilde
{\omega}_t}{\partial t}\>=\Delta_{\tilde{\omega}_t}(\frac{\partial u}{\partial t})-e^{-t}\<\tilde
{\omega}_t, \omega_0-\omega_\infty\>.$$
$$\frac{\partial}{\partial t}(\frac{\partial^2 u}{\partial t^2})=\<\tilde{\omega}_t, \frac{\partial^2
\tilde{\omega}_t}{\partial t^2}\>-(\frac{\partial\tilde{\omega}_t}{\partial t}, \frac{\partial\tilde
{\omega}_t}{\partial t})_{\tilde{\omega}_t}\leqslant\Delta_{\tilde{\omega}_t}(\frac{\partial
^2 u}{\partial t^2})+e^{-t}\<\tilde{\omega}_t, \omega_0-\omega_\infty\>.$$ 
Take summation of the above two and apply standard maximum principle argument to get: 
$$\frac{\partial}{\partial t}(\frac{\partial u}{\partial t}+u)\leqslant C,$$ 
and from this, it's easy to see 
$$\frac{\partial u}{\partial t}\leqslant C,$$ 
which gives the measure bound for ${\tilde{\omega}_t}^n=e^{\frac{\partial u}{\partial t}}
\Omega$. This allows us to apply the results of Kolodziej's (as in \cite{kojnotes} and 
\cite{koj06}) and our generalization (as in \cite{zh}) in respective situation of $[\omega_t]$, 
which provides the uniform bound for the metric potential along the flow after routine 
normalization in the cases under study. 

\begin{rema}

Even in the case of $[\omega_\infty]$ being K\"ahler, the result from pluripotential theory as 
in \cite{kojnotes} is used for the normalized metric potential bound, so the logic line of our 
argument is not quite the same as that in \cite{cao}. It would also be interesting to see whether 
the original argument of Cao's can also be carried through with fewer changes.  

\end{rema}  

We also need to derive some kind of low bound for $\frac{\partial u}{\partial t}$ in search of the 
possible metric bound (coming from volume and Laplacian controls). 

The following equation would be very useful later: 
\begin{equation}
\frac{\partial}{\partial t}(\frac{\partial u}{\partial t}+u)=\Delta_{\tilde{\omega}_t}(\frac{\partial u}
{\partial t}+u)-n+\frac{\partial u}{\partial t}+\<\tilde{\omega}_t, \omega_\infty\>.
\end{equation}   

\section{a Baby Version}

Let's start with the situation when there is no degeneration on the cohomology classes, i.e., 
$[\omega_\infty]$ is also K\"ahler. This is the first case in our main theorem, which could be 
seen as a natural generalization of Cao's work mentioned before after allowing the change 
of cohomology class along the flow. 

\subsection{Uniform Estimates and Global Existence}  

Global existence of the flow only needs estmates local in time as mentioned 
before. Now we have to go for estimates uniform for all time. The most essential 
part is the $C^0$ estimates. Up to now, there is only the lower bound for $\frac
{\partial u }{\partial t }$ left to be obtained. 

At first, let's assume $\omega_\infty>0$, which will not change the problem in any 
essential way. We'll remove this simplification later. \\      

As mentioned before, by the measure bound from the previous section, we already 
know that 
$$|v|\leqslant C$$ 
where $v=u-\int_Xu\Omega$. We have assumed $\int_X\Omega=1$ for simplicity 
of notation. 

We also know 
$$\frac{\partial v}{\partial t}=\frac{\partial u}{\partial t}-\int_X\frac{\partial u}{\partial t}
\Omega\geqslant \frac{\partial u}{\partial t}-C$$ 
from the upper bound of $\frac{\partial u}{\partial t}$, and so the lower bound of $\frac
{\partial u}{\partial t}$ would give that for $\frac{\partial v}{\partial t}$. 

In fact the inverse is also true as $(2.1)$ can be modified to be: 
$$\frac{\partial}{\partial t}(\frac{\partial u}{\partial t}+v)=\Delta_{\tilde{\omega}_t}(\frac
{\partial u}{\partial t}+v)-n+\frac{\partial v}{\partial t}+\<\tilde{\omega}_t, \omega_\infty
\>.$$

Assuming $\frac{\partial v}{\partial t}\geqslant -C$, we can get a lower bound for $\frac
{\partial u}{\partial t}$ by applying maximum principle as H. Tsuji did in \cite{tsu} using 
the control of volume by the control of trace, which is nothing but the classic 
algebraic-geometric mean value inequality. But we can actually do better by the following 
more delicate maximum principle argument. 

Consider the minimum value point, $p$, for $\frac{\partial u}{\partial t}+v$ for $X\times [0, T]$ 
with any fixed $0<T<\infty$. As usual, we only need to study the case when that point is not 
at the initial time since the situation for the initial time is well under control. At $p$, we have: 
$$n-\frac{\partial v}{\partial t}\geqslant \<\tilde{\omega}_t, \omega_\infty\>\geqslant n\cdot\bigl(
\frac{{\omega_\infty}^n}{{\tilde{\omega}_t}^n}\bigr)^{\frac{1}{n}}=n\cdot\bigl(\frac{{\omega_\infty}
^n}{e^{\frac{\partial u}{\partial t}}\Omega}\bigr)^{\frac{1}{n}}>0,$$
and so $(1-\frac{1}{n}\frac{\partial v}{\partial t})^n\cdot e^{\frac{\partial u}{\partial t}}\geqslant C
>0$. Using $\frac{\partial v}{\partial t}\geqslant\frac{\partial u}{\partial t}-C$, one arrives at: 
$$(C-\frac{\partial u}{\partial t})^n\cdot e^{\frac{\partial u}{\partial t}}\geqslant C>0$$
with $C-\frac{\partial u }{\partial t }>0$, which gives $\frac{\partial u}{\partial t}\geqslant -C$ at 
$p$. The uniform bound for $v$ thus gives the lower bound of $\frac{\partial u}{\partial t}$ over 
$X$ (and also for $\frac{\partial v}{\partial t}$) from the bound at $p$. \\   

There is another way of doing the maximum principle argument which might seems to be more 
direct in this case. This is also a very classic point of view when studying the flow. Basically, we 
examine the evolution of space-direction extremal value along the flow. This function, now only depending on time, would be (locally) Lipschitz simply by going through the definition, and so 
it'll be legitimate to consider the first order ordinary differential inequality. 

This kind of argument, in spirit, would be more delicate than what we used before. But actually 
for the differential inequality of interest here, the study would be as rough as before. Let's illutrate 
the idea below. \\ 

Set $A(t)={\rm min}_{X\times\{t\}}(\frac{\partial u}{\partial t}+v)$. Let's also take some $x(t)$ 
where the value $A(t)$ is achieved, but we do not assume (or need) any regularity of $x(t)$ 
with respect to $t$.  

Using the sign of Laplacian, we can derive the differential inequality for the function $A(t)$ 
as follows
\begin{equation}
\begin{split}
\frac{\partial A}{\partial t}
&\geqslant -n+\frac{\partial v}{\partial t}(x(t))+\<\tilde\omega_t, \omega_\infty\> (x(t))\\
&\geqslant -C+(\frac{\partial u}{\partial t}+v)(x(t))+Ce^{-(\frac{\partial u}{\partial t}+v)(x(t))}\\
&= -C+A+Ce^{-A}.\nonumber
\end{split} 
\end{equation} 

From this inequality, we can see that when $A$ is sufficiently small (i.e., very negative), 
$\frac{\partial A}{\partial t}$ would be big (i.e., very positive). It won't be hard to get a 
lower bound for $A$ from this mechanism. All the pieces from the previous argument 
are also used here, but this looks more straightforward (for people good at playing with 
ODE). 

\begin{rema}

Clearly, in the degenerate version of maximum principle as what will appear later, this 
point of view still works as long as the point $x(t)$ is in the regular part.

\end{rema}       

Then using classic second order estimate, we can have uniform control for the trace of 
$\tilde{\omega}_t$ (i.e., Lapacian). And so, together with the volume lower bound from 
the bound of $\frac{\partial u }{\partial t }$, we have them controlled uniformly as metric. 
Finally, high order derivatives are also uniformly controlled using classic estimates for 
parabolic PDE's (including Yau's computation and parabolic Schauder estimates). 
These are very standard arguments for the situation here. \\    

Now let's remove the assumption that $\omega_\infty >0$. We always have $\omega_
\infty+\sqrt{-1}\partial\bar{\partial}f>0$ for some smooth function $f$ over $X$ as $[\omega
_\infty]$ is K\"ahler. Also recall that $\omega_t=\omega_\infty+e^{-t}(\omega_0-\omega_
\infty)$. Let's now set 
$$\bar\omega_t=(\omega_\infty+\sqrt{-1}\partial\bar{\partial}f)+e^{-t}(\omega_0-(\omega_
\infty+\sqrt{-1}\partial\bar{\partial}f))=\omega_t+(1-e^{-t})\sqrt{-1}\partial\bar{\partial}f,$$ 
and clearly $\tilde\omega_t=\bar\omega_t+\sqrt{-1}\partial\bar{\partial}(u-(1-e^{-t})f)$. 
Define $w=u-(1-e^{-t})f$ and we have $\tilde \omega_t=\bar\omega_t+\sqrt{-1}\partial\bar
{\partial}w$. Clearly $\frac{\partial w}{\partial t}=\frac{\partial u}{\partial t}-e^{-t}f$ and taking 
$t$-derivative gives 
$$\frac{\partial}{\partial t}(\frac{\partial w}{\partial t})=\Delta_{\tilde\omega_t}(\frac{\partial w}
{\partial t})-e^{-t}\<\tilde\omega_t, \omega_0-\omega_\infty-\sqrt{-1}\partial\bar{\partial}f\>+
e^{-t}f.$$ 

Just as before, we need the following transformation of the above equation 
$$\frac{\partial}{\partial t}(\frac{\partial w}{\partial t}+\bar{w})=\Delta_{\tilde\omega_t}(\frac
{\partial w}{\partial t}+\bar{w})-n+\<\tilde\omega_t, \omega_\infty+\sqrt{-1}\partial\bar{\partial}
f\>+\frac{\partial\bar{w}}{\partial t}+e^{-t}f,$$
where $\bar{w}$ is the normalization of $w$ just as $v$ for $u$ before. 

We can now apply maximum principle for the above equaton at the (local in time) minimum 
value point of $\frac{\partial w}{\partial t}+\bar{w}$. At that point (if not at time $0$), we have 
$$n-\frac{\partial\bar{w}}{\partial t}\geqslant e^{-t}f+\<\tilde\omega_t, \omega_\infty+\sqrt
{-1}\partial\bar{\partial}f\>.$$

Without loss of generality, we can make sure $f>0$. So now one arrives at 
\begin{equation}
\begin{split}
n-\frac{\partial\bar{w}}{\partial t}
&\geqslant \<\tilde{\omega}_t, \omega_\infty+\sqrt{-1}\partial\bar{\partial}f\>\\
&\geqslant n\cdot\bigl(\frac{(\omega_\infty+\sqrt{-1}\partial\bar{\partial}f)^n}{{\tilde
{\omega}_t}^n}\bigr)^{\frac{1}{n}}=n\cdot\bigl(\frac{(\omega_\infty+\sqrt{-1}\partial
\bar{\partial}f)^n}{e^{\frac{\partial u}{\partial t}}\Omega}\bigr)^{\frac{1}{n}}>0, \nonumber
\end{split}
\end{equation}
which gives $(1-\frac{1}{n}\frac{\partial\bar{w}}{\partial t})^n\cdot e^{\frac{\partial u}
{\partial t}}\geqslant C>0$. We also have $\frac{\partial\bar{w}}{\partial t}\geqslant
\frac{\partial w}{\partial t}-C\geqslant\frac{\partial u}{\partial t}-C$, so we can have 
$$(C-\frac{\partial u}{\partial t})^n\cdot e^{\frac{\partial u}{\partial t}}\geqslant C>0$$
with $C-\frac{\partial u }{\partial t }>0$, and still conclude that $\frac{\partial u}{\partial t}
\geqslant -C$ at that point, and so $\frac{\partial w}{\partial t}\geqslant -C$ there.   

It's rather clear that $|\bar{w}|\leqslant C$ from the estimates for $v$ before since we 
do not need $\omega_\infty >0$ there yet. Hence we see $\frac{\partial w}{\partial t}+
\bar{w}\geqslant -C$ globally, which gives the uniform lower bound for $\frac{\partial 
w}{\partial t}$ and so for $\frac{\partial u}{\partial t}$ (as they differ only by a bounded 
term $e^{-t}f$). 

The argument for uniform higher derivatives is as standard as before.

\begin{rema}

The main philosophy of the above argument is that a choice of representative in a 
class boils down to terms like $f$ or $e^{-t}f$ for smooth function $f$ over $X$ which 
is clearly controlled along the flow and so should not bring any trouble. This observation 
is also useful when trying to apply Yau's Laplacian estimate (as in \cite{cao}) to get 
second order derivative control for the current situation.    

\end{rema}

Up to now, we have got the global existence of the flow and the uniformity of the estimates 
allows us to use the classic Ascoli-Azela's Theorem to get convergence for sequences of 
metrics along the flow. Just as in Cao's work, we should head for stronger convergence as 
discussed in next subsection.      

\subsection{Convergence} 

H. D. Cao's argument in \cite{cao} for convergence using Li-Yau's Harnack Inequality 
should be easy to get carried through here as for the equation:
$$\frac{\partial}{\partial t}(\frac{\partial u}{\partial t})
=\Delta_{\tilde{\omega}_t}(\frac{\partial u}{\partial t})-e^
{-t}\<\tilde{\omega}_t, \omega_0-\omega_\infty\>$$
since $\<\tilde{\omega}_t, \omega_0-\omega_\infty\>$ has been uniformly controlled, 
and so the extra term in comparison to Cao's situation is exponentially decreasing.  
Let's ilustrate some main points when adjusting his argument to the current situation 
in the following.\\  

For the exponential decreasing of the oscillation of $\frac{\partial u}{\partial t}$, we'll 
use Cao's argument for the following family of auxiliary functions: 
$$(\frac{\partial}{\partial t}-\Delta_{\tilde{\omega}_t})\phi_{T_0}=0, ~~\phi_{T_0}(T_0, 
\cdot)=\frac{\partial u}{\partial t}(T_0,\cdot)$$   
over $[T_0, \infty)\times X$ where $T_0\in [0, \infty)$. As we have already got the uniform 
estimates for $\frac{\partial u}{\partial t}$ and $\tilde{\omega}_t$, using Li-Yau's Harnack 
Inequality as Cao did, we have 
$${\rm osc}_X\phi_{T_0}(t)\leqslant Ce^{-a(t-T_0)}, ~~t\in[T_0, \infty)$$ 
where the positive constants are uniform for all $T$ as the metric control is uniform for 
all time.  
 
Using the uniform estimates along the flow as mentioned before, we have 
$$(\frac{\partial}{\partial t}-\Delta_{\tilde{\omega}_t})(\frac{\partial u}{\partial t}+Ce^{-t})
\leqslant 0,$$
$$(\frac{\partial}{\partial t}-\Delta_{\tilde{\omega}_t})(\frac{\partial u}{\partial t}-Ce^{-t})
\geqslant 0.$$

We also have the following equations
$$(\frac{\partial}{\partial t}-\Delta_{\tilde{\omega}_t})(\phi_{T-0}+Ce^{-T_0})=0,$$
$$(\frac{\partial}{\partial t}-\Delta_{\tilde{\omega}_t})(\phi_{T_0}-Ce^{-T_0})=0.$$

Comparing them and applying maximum principle, we get the decreasing of 
$${\rm max}_X(\frac{\partial u}{\partial t}+Ce^{-t}-\phi_{T_0}-Ce^{-T_0})$$ 
and the increasing of 
$${\rm min}_X(\frac{\partial u}{\partial t}-Ce^{-t}-\phi_{T_0}+Ce^{-T_0})$$ 
as time increases (starting from the time $T_0$). 

The values at $t=T_0$ for both quantities are $0$, so we have for $t\in [T_0, \infty)$, 
$$\frac{\partial u}{\partial t}\leqslant \phi_{T_0}+Ce^{-T_0}-Ce^{-t},$$
$$\frac{\partial u}{\partial t}\geqslant \phi_{T_0}-Ce^{-T_0}+Ce^{-t}.$$ 

Hence ${\rm osc}_X\frac{\partial u}{\partial t}\leqslant {\rm osc}_X\phi_{T_0}+
Ce^{-T_0}$ for $t\in [T_0, \infty)$. Using the result for $\phi_{T_0}$ stated above, 
we have ${\rm osc}_X\frac{\partial u}{\partial t}\leqslant Ce^{-a(t-T_0)}+Ce^{-T_0}$ 
for $t\geqslant T_0$. 
Taking $t=2T_0$ and noticing this is uniform for all $T_0$, we finally arrive at 
$${\rm osc}_X\frac{\partial u}{\partial t}\leqslant Ce^{-at}$$
for all time. Here the $a$ should differ from the previous one, but it's still a positive 
constant. 

This is exactly one of the essential results needed to draw the convergence for 
$t\to\infty$ as in Cao's.\\     

Set $\psi=\frac{\partial u }{\partial t }-\frac{\int_X\frac{\partial u }{\partial t }{\tilde
{\omega }_t }^n }{\int_X {\tilde{\omega }_t }^n}$. Clearly its difference from $\frac
{\partial v }{\partial t }$ is controlled by $Ce^{-at}$, but it is more convenient for the 
following consideration.  

We can have similar computation as in \cite{cao}, for the energy, 
$$E=\int_X\psi^2{\tilde{\omega }_t }^n,$$ 
to derive a differential inequality for it. There are more terms coming out than 
Cao's case, but they will all be terms controlled by $Ce^{-t}$ using the uniform 
estimates along the flow. Notice that though the volume is also changing 
along the flow, the variation is also well under control. In all, we get 
$$\frac{d E}{d t}\leqslant -CE+Ce^{-t}$$ 
for large $t$. The reason to get only for large $t$ is that we need the smallness 
of $\psi$ from the control of oscillator of $\frac{\partial u }{\partial t }$. From this 
differential inequality, we can still conclude the exponential decaying of $E$. 
\footnote{In fact, the exponential decaying of  $E$ can be deduced from the 
decaying of the oscillation of $\frac{\partial u}{\partial t}$ in a  more direct manner. 
But Cao's method above applying the differential inequality is more delicate and 
can easily be adjusted for higher order Sobolev estimates.} \\     

The final computation and argument of Cao to derive the $L^1$ convergence of 
the normalized metric potential can be carried through line by line in sight of the 
above results. Indeed, we can also justify the exponential convergence of the flow 
with little extra effort (just as what is in \cite{thesis}). 

\begin{rema}

In this situation, we now have a somewhat natural flow from one Ricci-flat metric to 
another Ricci-flat metric (in different K\"ahler classes of course) when $c_1(X)=
0$. Just need to choose $\Omega$ such that ${\rm Ric}(\Omega)=0$ for the flow.  

\end{rema}

\section{Main Interest: Degenerate Case}  

Of course, our main interest is when $[\omega_\infty]$ is degenerate as K\"ahler 
class. In \cite{thesis}, we have discussed the corresponding Monge-Amp\`ere 
equation using other perturbation for method of continuity. Now we want to see 
whether the modified flow can help us to construct a solution for the Monge-Amp\`ere 
equation (as the limiting equation). At this moment, our manifold $X$ is assumed 
to be projective to get into algebraic geometry context for the notions of semi-ample 
and big. 

As discussed before, we have the existence of the smooth flow as long as $[\omega_t]$ 
remains K\"ahler. There are two cases, i.e., up to infinite time and up to finite time. We 
discuss them separately and finish the proof of the theorem.   

\subsection{Infinite Time Case} 

Let's assume here that $[\omega_\infty]$ is semi-ample and big. We still have the $L^\infty$ 
bound of the normalized metric potential $v$ as before using the result on degenerate 
Monge-Amp\`ere equation (from \cite{ey-gu-ze} and \cite{zh}). Now $(2.1)$ can be modified 
as: 
$$\frac{\partial}{\partial t}(\frac{\partial u}{\partial t}+v-\epsilon{\rm log}|\sigma|^2)=\Delta_
{\tilde{\omega}_t}(\frac{\partial u}{\partial t}+v-\epsilon{\rm log}|\sigma|^2)-n+\frac{\partial v}
{\partial t}+\<\tilde{\omega}_t, \omega_\infty+\epsilon\sqrt{-1}\partial\bar{\partial}{\rm log}
|\sigma|^2\>$$
with $\omega_\infty+\epsilon\sqrt{-1}\partial\bar{\partial}{\rm log}|\sigma|^2>0$, where the 
positive $\epsilon$ can be as close to $0$ as possible.  The introduction of a singular term 
like this, as far as I know, was initiated by Tsuji in \cite{tsu}, which gives a natural and simple 
description of an algebraic geometry fact in analysis of the related PDE's. 

The following classic results from algebraic geometry are useful for us. See in \cite{kawa1} 
and \cite{kawa3} for related discussion. The second one is called Kodaira's Lemma as in 
\cite{tsu2} and will be applied for the finite time case later discussed. The point of view for 
translating these results to the analytic statement as above is very standard as described in 
\cite{gri-har}.  

\begin{lemm}

Let $L$ be a divisor in a projective manifold $X$. If $L$ is nef. and big, then there is an 
effective divisor $E$ and a number $a>0$ such that $L-\epsilon E$ is K\"ahler for any
$\epsilon \in (0,a)$.

\end{lemm}

\begin{lemm}

Let $L$ be a divisor in a projective manifold $X$. If $L$ is big, then there is an effective 
divisor $E$ such that $L-\epsilon E$ is K\"ahler for $\epsilon \in (a,b)$ where $0\leqslant 
a<b<\infty$.

\end{lemm}
 
Similar argument as before would give a degenerate lower bound \footnote{The positive 
constant $C$ below might depend on the other positive constant $\epsilon$. Hopefully, 
this won't bring any confusion.} as 
$$\frac{\partial u}{\partial t}\geqslant -C+\epsilon{\rm log}|\sigma|^2.$$ 
Basically, we still have $\frac{\partial v }{\partial t }\geqslant\frac{\partial u }{\partial t }-C$. 
Then considering the minimum value point of the term naturally considered by the equation 
above, we know $\frac{\partial u }{\partial t}$ could not be too small at that point using the 
contradiction as for the baby version, which would essentially give the bound claimed above.  

Now the degenerate second order estimate and high ones would still be OK by the standard 
procedure. More specifically, for the second order estimate, one considers the following 
equation
$$(\tilde{\omega}_t+\epsilon\sqrt{-1}\partial\bar{\partial }{\rm log }|\sigma |^2+\sqrt{-1}
\partial\bar{\partial }(v-\epsilon\,{\rm log }|\sigma |^2))^n=e^{\frac{\partial u}{\partial t} }
\Omega.$$ 

Applying Yau's computation in \cite{yau} and using degenerate maximum principle 
argument as in \cite{t-znote}, we can get the degenerate Laplacian bound. \footnote
{In fact, one only needs the uniform upper bound of $\frac{\partial u}{\partial t}$ to carry 
through the Laplacian estimate by noticing the dominance of $e^{-\frac{1}{n}\frac
{\partial u}{\partial t}}$ over $-\frac{\partial u}{\partial t}$ when $\frac{\partial u}{\partial t}$ 
is small.} Combining with the degenerate control for volume, we have achieved local 
(or degenerate) bound for metrics along the flow. 

The treatment for higher derivatives would be standard. We provide some details at the 
last section.   

\begin{rema}

There is a big difference from the situation in \cite{t-znote} which we want to point out. 
The metric potential along the flow can be bounded (though in a degenerate way) simply 
from the flow argument, but we can not do that here at this moment. The bound for 
(normalized) metric potential is coming from results proved by arguments in pluripotential 
theory. That's why we need semi-ample (not just nef.) here.  

\end{rema}

Though our estimates are uniform for all time now, which gives sequence convergence 
for the flow, there is still this big issue about convergence along the flow which is crucial 
to describe the limit itself. As discussed in the baby version, the counterpart in \cite{cao} 
makes use of Li-Yau's Harnack Inequality, which can be applied for the non-degenerate 
case as before. But the situation right now is very different. It seems to me that new 
method needs to be introduced for this purpose. Let's make the conjecture about the flow 
convergence. 

\begin{conj}

For $[\omega_\infty]$ semi-ample and big, as $t\to \infty$, this modified K\"ahler-Ricci flow 
converges weakly over $X$ and locally smoothly out of the stable base locus set of this 
cohomology class to the unique (bounded) solution of the limiting degenerate 
Monge-Amp\`ere equation. 

\end{conj}

We can prove that for infinite time case, the volume form has uniform lower bound for 
all time as stated in Theorem 1.2. This might help to get the convergence of the flow 
and is also a nice application of a similar result for the following more canonical 
K\"ahler-Ricci flow. 

Set $\hat{\omega}_t=\omega_t+\sqrt{-1}\partial\bar{\partial}\phi$. In the level of potential, 
consider the flow
$$\frac{\partial\phi}{\partial t}={\rm log}\frac{{\hat{\omega}_t}^n}{\Omega}-\phi, ~~~~\phi (0, 
\cdot)=0.$$
The corresponding flow in the level of metric is the following, 
$$\frac{\partial\hat{\omega}_t}{\partial t}=-{\rm Ric}(\hat
{\omega}_t)+{\rm Ric}(\Omega)-\hat{\omega}_t+\omega_\infty, 
~~~~\hat{\omega}_0=\omega_0.$$

In the case, we have $[\omega_\infty]$ is big and semi-ample, so as discussed in 
\cite{t-znote} and \cite{thesis}, the following controls are available, \footnote{The 
lower bound of $\frac{\partial\phi}{\partial t}$ is in \cite{thesis}. Basically one makes 
use of the essential decreasing of the volume form and the fact that for infinite time 
limit, the derivative of potential has to go to $0$ (in the regular part). }
$$|\phi|\leqslant C, ~~|\frac{\partial \phi}{\partial t}|\leqslant C,$$  
which give a lower bound for the volume form $\hat{\omega}_t$ for all time and 
we are looking for a similar thing for ${\tilde{\omega}_t}^n$.

\begin{rema}

The uniform volume lower bound is a pretty interesting fact as the class $[\omega_\infty]$ 
is not K\"ahler, but somehow we have that $[\omega_\infty]^n>0$ also makes sense in a 
pointwise fashion. 

\end{rema}

Let's recall the following equations used before for the flow 
considered in this note. $v$ is the normalization of $u$ as before.  
$$\frac{\partial}{\partial t}(\frac{\partial u}{\partial t})=\Delta(\frac
{\partial u}{\partial t})-\<\tilde{\omega}_t, e^{-t}(\omega_0-\omega
_\infty)\>,$$
$$\frac{\partial}{\partial t}(\frac{\partial u}{\partial t})=\Delta(\frac
{\partial u}{\partial t}+v)-n+\<\tilde{\omega}_t, \omega_\infty\>.$$

Fix some constant $T_1>0$, product the first equation with $e^{-T_1}$ 
and taking the difference of them, we have 
$$\frac{\partial}{\partial t}((1-e^{-T_1})\frac{\partial u}{\partial t})=
\Delta\bigl((1-e^{-T_1})\frac{\partial u}{\partial t}+v\bigr)-n+\<\tilde
{\omega}_t, \omega_{t+T_1}\>.$$

Now using the solution for the other flow, $\phi$, this equation can 
be transformed as follows 
$$\frac{\partial}{\partial t}((1-e^{-T_1})\frac{\partial u}{\partial t})=
\Delta\bigl((1-e^{-T_1})\frac{\partial u}{\partial t}+v-\phi(t+T_1)\bigr)
-n+\<\tilde{\omega}_t, \hat{\omega}_{t+T_1}\>$$
with some emphasize on the time parameter. The Laplacian 
is still with respect to the metric $\tilde{\omega}_t$. In similar 
spirit as before, we modify the equation to be 
\begin{equation}
\begin{split}
\frac{\partial}{\partial t}((1-e^{-T_1})\frac{\partial u}{\partial t}+v-
\phi(t+T_1))
&= \Delta((1-e^{-T_1})\frac{\partial u}{\partial t}+v-\phi(t+T_1))-n+
\frac{\partial v}{\partial t} \\
&~~~~-\frac{\partial\phi(t+T_1)}{\partial t}+\<\tilde{\omega}_t, \hat
{\omega}_{t+T_1}\>. \nonumber
\end{split}
\end{equation}

Let $A=(1-e^{-T_1})\frac{\partial u}{\partial t}+v-\phi(t+T_1)$ and using 
the following known estimates
$$\frac{\partial v}{\partial t}\geqslant \frac{\partial u}{\partial t}-C, 
~~\frac{\partial\phi(t+T_1)}{\partial t}\geqslant -C, ~~{\hat{\omega}
_t}^n\geqslant C\Omega,$$ 
one arrives at 
$$\frac{\partial A}{\partial t}\geqslant \Delta A+\frac{\partial u}
{\partial t}-C+C\cdot e^{-\frac{1}{n}\frac{\partial u}{\partial t}}.$$

Use similar maximum principle argument as before, one can 
conclude the lower bound for $A$, and so for $\frac{\partial u}
{\partial t}$, which gives the lower bound for the volume form 
${\tilde{\omega}_t}^n$. 

\begin{rema}

The translation of time by $T_1$ makes the infinite time situation 
special. In comparison, we do not have uniform volume lower 
bound for finite time case for both flows (at least at this moment). 

\end{rema}

\subsection{Finite Time Limit} 

Now we consider the case when $[\omega_\infty]$ is only big. More specifically, 
the flow only exists up to some finite time $T$. We also require $[\omega_T]$, 
which is clearly nef. and big,  to be semi-ample. \footnote{This is not such a 
horrible assumption as it is the case, when $[\omega_\infty]=K_X$ and $[\omega
_0]$ is rational, from algebraic geometry results.} 

Let's first consider the situation roughly. Those degenerate estimates would still 
be available, though the $\epsilon$ can't be too small now in sight of Lemma 4.2. 
The advantage about finite time is that the metric potential $u$ is (degenerately) 
bounded by itself (without normalization) using the bound for its time derivative, 
and it'll also be decreasing after controllable normalization (by $-Ct$) in sight of 
the uniform upper bound for $\frac{\partial u }{\partial t }$. So as in \cite{t-znote}, 
the (local) convergence for $t\to T$ is achieved. 

This local convergence would be out of the "stable base locus set" of $[\omega_
\infty]$. Clearly, it would be more satisfying to get this with respect to $[\omega_T]$. 
We can do this in the same way in which we can also improve the result in 
\cite{t-znote}. Simply speaking, we can use a virtual time. Let's get the crucial 
estimate for $\frac{\partial u }{\partial t }$ below.  

We can easily have the following two equations. 
$$\frac{\partial}{\partial t}(\frac{\partial u}{\partial t}+v)=\Delta_{\tilde{\omega}_t}
(\frac{\partial u}{\partial t}+v)-n+\frac{\partial v}{\partial t}+\<\tilde{\omega}_t, 
\omega_\infty\>,$$
$$\frac{\partial}{\partial t}(e^{t-T}\frac{\partial u}{\partial t})=\Delta_{\tilde{\omega}_
t}(e^{t-T}\frac{\partial u}{\partial t})+e^{t-T}\frac{\partial u }{\partial t }-e^{-T}\<\tilde
{\omega}_t, \omega_0-\omega_\infty\>.$$

Take difference to get
$$\frac{\partial}{\partial t}((1-e^{t-T})\frac{\partial u}{\partial t}+v)=\Delta_{\tilde
{\omega}_t}((1-e^{t-T})\frac{\partial u}{\partial t}+v)-n+\frac{\partial v}{\partial t}-
e^{t-T}\frac{\partial u }{\partial t }+\<\tilde{\omega}_t, \omega_T\>.$$

As before, take $\sigma$ for $[\omega_T]$ such that $\omega_T+\epsilon
\sqrt{-1}\partial\bar{\partial }{\rm log }|\sigma |^2>0$ for any positive $\epsilon$ 
small enough. Using it to perturb the above equation, one arrives at
\begin{equation}
\begin{split}
&~~ \frac{\partial}{\partial t}((1-e^{t-T})\frac{\partial u}{\partial t}+v-\epsilon\,
{\rm log }|\sigma |^2)\\
&= \Delta_{\tilde{\omega}_t}((1-e^{t-T})\frac{\partial u}{\partial t}+v-\epsilon\,
{\rm log }|\sigma |^2)-n+\frac{\partial v}{\partial t}-e^{t-T}\frac{\partial u }
{\partial t }+\<\tilde{\omega}_t, \omega_T+\epsilon\sqrt{-1}\partial\bar{\partial }{\rm log }
|\sigma |^2\>.\nonumber
\end{split}
\end{equation}

Set $A=(1-e^{t-T})\frac{\partial u}{\partial t}+v-\epsilon\,{\rm log }|\sigma |^2$. 
We can have, out of $\{\sigma=0\}$,   
$$\frac{\partial A }{\partial t }\geqslant \Delta_{\tilde{\omega }_t}A-C+(1-e^{t-T})
\frac{\partial u }{\partial t }+Ce^{-\frac{1}{n}\frac{\partial u }{\partial t } }.$$
Then let's do the maximum principle argument. Recall that the time $t\in[0, T)$.
At the minimum value point of $A$ (assuming it is not at the initial time), which 
is clearly out of $\{\sigma=0\}$, we can see $\frac{\partial u}{\partial t }$ can not 
be too small (negative). Thus $A$ can not be too small there, either. That gives
$$(1-e^{t-T})\frac{\partial u}{\partial t}+v-\epsilon\,{\rm log }|\sigma |^2\geqslant 
-C.$$

The problem coming from the fact that $1-e^{t-T}$ would go to $0$ as $t\to T$ 
can be solved by using a "virtual" time $T_\epsilon>T$ which satisfies $\omega_
{T_\epsilon}+\epsilon\sqrt{-1}\partial\bar{\partial }{\rm log }|\sigma |^2>0$ for 
some fixed $\epsilon>0$.  Then the same estimate
$$(1-e^{t-T_\epsilon })\frac{\partial u}{\partial t}+v-\epsilon\,{\rm log }|\sigma |^2
\geqslant -C$$
which gives
$$\frac{\partial u }{\partial t }\geqslant -C_\epsilon+C_\epsilon{\rm log }|\sigma 
|^2.$$

Notice that now the $\sigma$ is for the class $[\omega_T]$ and so we can conclude 
the local convergence out of the stable base locus set of $[\omega_T]$. \\

One might also want to do the maximum principle argument in another flavor 
just as what is done for the baby version. We have to do it more carefully as follows. 

The differential inequality for $A$ is  
$$\frac{\partial A }{\partial t }\geqslant \Delta_{\tilde{\omega }_t}A-C+(1-e^{t-T})\frac
{\partial u }{\partial t }+Ce^{-\frac{1}{n}\frac{\partial u }{\partial t } }.$$

One wants to change the last two terms to functions on $A$ with the right direction 
of control. 

The last term can be treated with ease as $e^{t-T}\frac{\partial u }{\partial t}-v+\epsilon
\,{\rm log }|\sigma |^2\leqslant C$, but it won't be so easy for the other term as $-\epsilon
\,{\rm log }|\sigma |^2$ can not be bound from above over $X$ by any constant. In fact, 
the trick is to treat them together. 

Set $B=(1-e^{t-T})\frac{\partial u}{\partial t}+v$ and we have  
$$(1-e^{t-T})\frac{\partial u}{\partial t}+Ce^{-\frac{1}{n}\frac{\partial u}{\partial t}}\geqslant 
-C+B+Ce^{-\frac{B}{n}}.$$

The function (over $B$), $B+Ce^{-B}$ would be decreasing with respect to $B$ for 
small enough $B$ by derivative consideration. And so for $B$ small enough (i.e., 
$(1-e^{t-T})\frac{\partial u }{\partial t }$ small enough), we can change $B$ to $A=
B-\epsilon\,{\rm log }|\sigma |^2$ and that should do it. \\  

The proof of Theorem 1.2 is finished. 

\section{Higher Order Estimates}

We provide a short discussion on the degenerate third and higher estimates for 
K\"ahler-Ricci flow over a closed (algebraic) manifold, $X$. It works for the modified 
flow here and others (as the one in \cite{t-znote} appearing in Subsection 4.2) 
with $[\omega_\infty]$ being big.  

The flow equation on the potential level is 
$$\frac{\partial u}{\partial t}={\rm log}\frac{(\omega_t+\sqrt{-1}\partial\bar{\partial}u)^n}
{\Omega}-u, ~~~~u(0, \cdot)=0,$$
or without the $-u$ term on the right hand side, where $\omega_t=\omega_\infty+e^{-t}
(\omega_0-\omega_\infty)$ with $\omega_0$ being the initial K\"ahler metric, $\omega
_\infty$ being a smooth representative for the (formal) infinite limiting class and $\Omega$ 
being a smooth volume for over . $\tilde\omega_t=\omega_t+\sqrt{-1}\partial\bar{\partial}
u$ is the metric solution of the flow.  

The class $[\omega_\infty]$ is big and $T\leqslant\infty$ is the singular time from 
cohomology concern with $[\omega_T]$ being nef. and big. The following estimates 
are available: 
$$\epsilon{\rm log}|\sigma|^2-C_\epsilon \leqslant u\leqslant C, ~~C{\rm log}
|\sigma|^2-C\leqslant \frac{\partial u}{\partial t}\leqslant C, ~~\<\omega_0, \tilde
\omega_t\>\leqslant C|\sigma|^{-l},$$
where $E=\{\sigma=0\}$ is a proper chosen devisor such that $[\omega_T]-\epsilon 
E$ is K\"ahler. $\sigma$ is a holomorphic section of the line bundle, and so with a 
fixed hermitian metric, $|\sigma|^2$ is a smooth function valued in $[0, C]$. 
 
The higher estimates are discussed briefly to achieve the full local regularity. 
Here we would like to go for the third order estimate a little more carefully. Then 
the rest follows from parabolic version of Schauder estimates in a standard way.
Yau's computation in \cite{yau} is what we need. 

As in Yau's computation, the term $S=\tilde g^{i\bar j}\tilde g^{k\bar l}
\tilde g^{m\bar n}u_{i\bar l m}u_{\bar j k\bar n}$ is considered, where 
the covariant derivative is with respect to uniform "background" 
metric. 

If the flow metric control is uniform, then the parabolic version of Yau's computation is
$$(\Delta_{\tilde\omega_t}-\frac{\partial}{\partial t})S\geqslant -C\cdot S
-C.$$

To adjust the result, one only need to see the metric is controlled (uniformly in time) 
as follows 
$$|\sigma|^\beta\omega_0\leqslant \tilde\omega_t\leqslant |\sigma|^{-\beta}\omega_0$$ 
for large positive constant $\beta$. 

Then we know by very carefully going through Yau's computation that
$$|\sigma|^{2N}(\Delta_{\tilde\omega_t}-\frac{\partial}{\partial t})S\geqslant -C|\sigma|^
{2N-\beta}\cdot S-C$$
with $N$ chosen large enough to dominate all degenerate terms. \\

Of course, now we want to see how $|\sigma|^{2N} S$ is acted by the 
heat operator. The only additional part is from the action of  $\Delta_{
\tilde\omega_t}$. There are two terms. One is clearly $2{\rm Re}(\nabla 
|\sigma|^{2N}, \nabla S)_{\tilde\omega_t}$. The other one is $\Delta_{
\tilde\omega_t}|\sigma|^{2N}\cdot S$. 

For the first one, 
$\nabla S=\nabla (|\sigma|^{2N}S\cdot |\sigma|^{-2N})=|\sigma|^{-2N}
\nabla (|\sigma|^{2N}S)-N|\sigma|^{-2}S\nabla |\sigma|^2$.

For the second one,
\begin{equation}
\begin{split} 
\Delta_{\tilde\omega_t}|\sigma|^{2N}
&=  \Delta_{\tilde\omega_t}(e^{N{\rm log}|\sigma|^2})=(|\sigma|^{2N}
N({\rm log}|\sigma|^2)_{\bar i})_i \\
&= N^2|\sigma|^{2N}|\nabla{\rm log}|\sigma|^2|^2+N|\sigma|^{2N}
\<\tilde\omega_t, \sqrt{-1}\partial\bar\partial{\rm log}|\sigma|^2\> \\
&\geqslant -N|\sigma|^{2N}\<\tilde\omega_t, -\sqrt{-1}\partial\bar\partial
{\rm log}|\sigma|^2\> \nonumber
\end{split}
\end{equation}

Out of $E=\{\sigma=0\}$, $-\sqrt{-1}\partial\bar\partial{\rm log}|\sigma|^2$ 
is nothing but the curvature form of the corresponding line bundle, still 
denoted by $E$. Using the degenerate metric bound, one has 
$$\Delta_{\tilde\omega_t} |\sigma|^{2N} \geqslant -N|\sigma|^{2N}
\<\tilde\omega_t, E\>\geqslant -C|\sigma|^{2N-\beta},$$

\begin{rema}

If we are in semi-ample case, with proper choice of the hermitian metric 
for the bundle ($|\cdot|$ above), we can make sure that $\Phi-\epsilon 
E>0$ (since the corresponding cohomology class is K\"ahler) where 
$\Phi$ is the pullback of a K\"ahler metric from the image of the map 
which is constructed from the semi-ample class $[\omega_T]$. In this 
case, we also have better zero order bounds. Moreover, using Schwarz 
type of estimates as in \cite{scalar-curvature}, we can have $\<\tilde
\omega_t, \Phi\>\leqslant C$. So now 
$$\Delta_{\tilde\omega_t} |\sigma|^{2N} \geqslant -N|\sigma|^{2N}
\<\tilde\omega_t, E\>\geqslant -C|\sigma|^{2N}\<\tilde\omega_t, \Phi\>
\geqslant -C|\sigma|^{2N},$$
which is better here but not going to make too much difference as one 
continues. 

\end{rema}

Anyway, we arrive at 
\begin{equation}
\begin{split}
(\Delta_{\tilde\omega_t}-\frac{\partial}{\partial t})(|\sigma|^{2N}S)
&\geqslant -C|\sigma|^{2N-\beta}\cdot S-C \\
&~~~ +2{\rm Re}\bigl(\nabla |\sigma|^{2N}, |\sigma|^{-2N}\nabla 
(|\sigma|^{2N}S)-N|\sigma|^{-2}S\nabla |\sigma|^2\bigr)_{\tilde
\omega_t} \\
&= -C|\sigma|^{2N-\beta}S-C+2{\rm Re}\bigl(\nabla ({\rm log}
|\sigma|^{2N}), \nabla(|\sigma|^{2N}S)\bigr)_{\tilde\omega_t} \\
&~~~ -N^2|\sigma|^{2N-4}S|\nabla |\sigma|^2|^2. \\
&\geqslant -C|\sigma|^{2N-2-\beta}S-C+2{\rm Re}\bigl(\nabla 
({\rm log}|\sigma|^{2N}), \nabla(|\sigma|^{2N}S)\bigr)_{\tilde
\omega_t}, \nonumber
\end{split}
\end{equation}
where we use $|\nabla|\sigma|^2|^2\leqslant C|\sigma|^{2-\beta}$ 
for the last step \footnote{The outside $|\cdot|$ is $\tilde\omega_t$.}. 

Also as in \cite{yau}, we consider the $\<\omega_{t, \epsilon}, \tilde\omega_t\>$ 
acted by the heat operator where $\omega_{t, \epsilon}$ is the perturbation for 
the "background" form. 

Had the metric control been uniform, one has
$$(\Delta_{\tilde\omega_t}-\frac{\partial}{\partial t})(\<\omega_{t, 
\epsilon}, \tilde\omega_t\>)\geqslant C\cdot S-C.$$ 

For our case, similar to $S$, we have instead 
$$|\sigma|^{2N}(\Delta_{\tilde\omega_t}-\frac{\partial}{\partial t})
(\<\omega_{t, \epsilon}, \tilde\omega_t\>)\geqslant C|\sigma|^{2N
+\beta}S-C.$$ 

The exact same procedure as done for $S$ above gives us
\begin{equation}
\begin{split}
&~~ (\Delta_{\tilde\omega_t}-\frac{\partial}{\partial t})(|\sigma|^{2N}
\<\omega_{t, \epsilon}, \tilde\omega_t\>) \\
&\geqslant C|\sigma|^{2N+\beta}\cdot S-C-C|\sigma|^{2N-\beta}
\<\omega_{t, \epsilon}, \tilde\omega_t\> \\
&~~~ +2{\rm Re}\bigl(\nabla |\sigma|^{2N}, |\sigma|^{-2N}\nabla 
(|\sigma|^{2N}\<\omega_{t, \epsilon}, \tilde\omega_t\>)-N|\sigma|^
{-2}\<\omega_{t, \epsilon}, \tilde\omega_t\>\nabla |\sigma|^2\bigr)_
{\tilde\omega_t} \\
&\geqslant C|\sigma|^{2N+\beta}S-C+2{\rm Re}\bigl(\nabla ({\rm log}
|\sigma|^{2N}), \nabla(|\sigma|^{2N}\<\omega_{t, \epsilon}, \tilde
\omega_t\>)\bigr)_{\tilde\omega_t} \\
&~~~ -C|\sigma|^{2N-2-\beta}\<\omega_{t, \epsilon}, \tilde\omega_t\>. 
\nonumber
\end{split}
\end{equation}

Properly choosing large constants $N_1>N_2>0$ and $C$'s, we have  
\begin{equation}
\begin{split}
&~~ (\Delta_{\tilde\omega_t}-\frac{\partial}{\partial t})(|\sigma|^{2N_1}S
+C|\sigma|^{2N_2}\<\omega_{t, \epsilon}, \tilde\omega_t\>) \\
&\geqslant C|\sigma|^{2N_2+\beta}\cdot S-C-C|\sigma|^{2N_2-2-\beta}
\<\omega_{t, \epsilon}, \tilde\omega_t\> \\
&~~~ +2{\rm Re}\bigl(\nabla ({\rm log}|\sigma|^{2N_1}), \nabla(|\sigma|^
{2N_1}S)\bigr)_{\tilde\omega_t}+2{\rm Re}\bigl(\nabla ({\rm log}|\sigma|^
{2N_2}), \nabla (C|\sigma|^{2N_2}\<\omega_{t, \epsilon}, \tilde\omega_t\>)
\bigr)_{\tilde\omega_t}. \nonumber
\end{split}
\end{equation}

For $N_1$ and $N_2$, only need $2N_1-2-\beta\geqslant 2N_2+\beta$ 
at this moment.  But they will be fixed later and large. \\ 

Now apply maximum principle argument. At the (local in time) maximum 
point,  for $|\sigma|^{2N_1}S+C|\sigma|^{2N_2}\<\omega_{t, \epsilon}, \tilde
\omega_t\>$, which clearly exists out of $\{\sigma=0\}$ and assume is not 
at the initial time, one has $\nabla(|\sigma|^{2N_1}S)=-\nabla(C|\sigma|^
{2N_2}\<\omega_{t, \epsilon}, \tilde\omega_t\>)$ and
\begin{equation}
\begin{split}
0
&\geqslant C|\sigma|^{2N_2+\beta}\cdot S-C-C|\sigma|^{2N_2-2-\beta}
\<\omega_{t, \epsilon}, \tilde\omega_t\> \\
&~~~ +2{\rm Re}\bigl(\nabla ({\rm log}|\sigma|^{2N_1}), \nabla(|\sigma|^
{2N_1}S)\bigr)_{\tilde\omega_t}+2{\rm Re}\bigl(\nabla ({\rm log}|\sigma|^
{2N_2}), \nabla (C|\sigma|^{2N_2}\<\omega_{t, \epsilon}, \tilde\omega_t\>)
\bigr)_{\tilde\omega_t} \\
&= C|\sigma|^{2N_2+\beta}\cdot S-C-C|\sigma|^{2N_2-2-\beta}\<\omega_
{t, \epsilon}, \tilde\omega_t\> \\
&~~~ +2{\rm Re}\bigl(-\nabla ({\rm log}|\sigma|^{2N_1})+\nabla ({\rm log}
|\sigma|^{2N_2}), \nabla (C|\sigma|^{2N_2}\<\omega_{t, \epsilon}, \tilde
\omega_t\>)\bigr)_{\tilde\omega_t} \\
&\geqslant C|\sigma|^{2N_2+\beta}\cdot S-C-C|\bigl(\nabla ({\rm log}
|\sigma|^{2}), \nabla (|\sigma|^{2N_2}\<\omega_{t, \epsilon}, \tilde
\omega_t\>)\bigr)_{\tilde\omega_t}|. \nonumber
\end{split}
\end{equation}

For the last term, we have
\begin{equation}
\begin{split}
&~~ |\bigl(\nabla ({\rm log}|\sigma|^{2}), \nabla (|\sigma|^{2N_2}\<\omega_
{t, \epsilon}, \tilde\omega_t\>)\bigr)_{\tilde\omega_t}| \\
&= |\bigl( |\sigma|^{-2}\nabla |\sigma|^2, N_2|\sigma|^{2N_2-2}\nabla 
|\sigma|^2\<\omega_{t, \epsilon}, \tilde\omega_t\>+|\sigma|^{2N_2}\nabla
\<\omega_{t, \epsilon}, \tilde\omega_t\>\bigr)_{\tilde\omega_t}| \\
&\leqslant   |\bigl( |\sigma|^{-2}\nabla |\sigma|^2, N_2|\sigma|^{2N_2-2}
\nabla |\sigma|^2\<\omega_{t, \epsilon}, \tilde\omega_t\>\bigr)_{\tilde
\omega_t}|+ |\bigl( |\sigma|^{-2}\nabla |\sigma|^2, |\sigma|^{2N_2}\nabla
\<\omega_{t, \epsilon}, \tilde\omega_t\>\bigr)_{\tilde\omega_t}| \\
&\leqslant C|\sigma|^{2N_2-2-2+1+1-\beta-\beta}+|\sigma|^{2N_2-2}|\nabla 
|\sigma|^2|\cdot |\nabla\<\omega_{t, \epsilon}, \tilde\omega_t\>| \\ 
&\leqslant C|\sigma|^{2N_2-2-2\beta}+C|\sigma|^{2N_2-2+1-\frac{\beta}{2}}
\cdot |\nabla\<\omega_{t, \epsilon}, \tilde\omega_t\>|
\nonumber
\end{split}
\end{equation}

Now one needs to realize that 
$$|\nabla\<\omega_{t, \epsilon}, \tilde\omega_t\>|=|\nabla (F+\Delta_
{\omega_{t,\epsilon}}u)|\leqslant |\nabla F|+|\nabla \Delta_{\omega_
{t,\epsilon}}(u)|\leqslant |\sigma|^{-\frac{\beta}{2}}+C|\sigma|^{-2\beta}
S^{\frac{1}{2}}$$
with $F$ being a well controlled function. 

Combining all this, we have at that maximum point, 
$$0\geqslant C|\sigma|^{2N_2+\beta}\cdot S-C-C|\sigma|^{2N_2-2-2\beta}
-C|\sigma|^{2N_2-1-\beta}-C|\sigma|^{2N_2-1-\frac{5\beta}{2}}\cdot S^{\frac
{1}{2}}.$$
For large enough $N_2$, we have,
$$0\geqslant |\sigma|^{2N_2+\beta}\cdot S-C (|\sigma|^{2N_2
+\beta}\cdot S)^{\frac{1}{2}}-C,$$
and so $|\sigma|^{2N_2+\beta}\cdot S\leqslant C$.
For $N_1$ even larger, we have uniform upper bound for $|\sigma|^{2N_1}
S+C|\sigma|^{2N_2}\<\omega_{t, \epsilon}, \tilde\omega_t\>$ at that point 
and so it is true globally which provides the bound
$$S\leqslant C|\sigma|^{-2N_1}.$$

This gives local $C^{2, \alpha}$ bound for the metric along the flow, then 
parabolic version of Schauder estimates carry though to provide all the 
local higher order bounds.


\begin{thebibliography}{$$}

\bibitem{cao} Cao, Huaidong: Deformation of Kaehler 
metrics to Kaehler-Einstein metrics on compact Kaehler 
manifolds. Invent. Math. 81(1985), no. 2, 359--372.

\bibitem{ey-gu-ze} Philippe Eyssidieux; Vincent Guedj; 
Ahmed Zeriahi: Singular K\"ahler-Einstein metrics. 
ArXiv, math/0603431.  

\bibitem{gri-har} Griffiths, Phillip; Harris, Joseph: 
Principles of algebraic geometry. Pure and Applied 
Mathematics. Wiley-Interscience [John Wiley \& Sons], 
New York, 1978. xii+813pp.

\bibitem{ham} Hamilton, Richard S.: Three-manifolds with 
positive Ricci curvature. J. Differential Geom. 17 (1982), no. 
2, 255-306. 

\bibitem{kawa1} Kawamata, Yujiro: The cone of curves 
of algebraic varieties. Ann. of Math. (2) 119 (1984), no.3, 
603--633.

\bibitem{kawa3} Kawamata, Yujiro: A generalization 
of Kodaira-Ramanujam's vanishing theorem. Math. Ann. 
261 (1982), no. 1, 43--46.

\bibitem{kojnotes} Kolodziej, Slawomir: The complex 
Monge-Amp\`ere equation and pluripotential theory. Mem. Amer. 
Math. Soc. 178 (2005), no. 840, x+64 pp. 

\bibitem{koj06} Kolodziej, Slawomir: H\"older continuity 
of solutions to the complex Monge-Amp\`ere equation with the right 
hand side in $L^p$. Preprint. 

\bibitem{t-znote} Tian, Gang; Zhang, Zhou: On the 
K\"ahler-Ricci flow on projective manifolds of general 
type. Chinese Annals of Mathematics - Series B, Volume 27, 
Number 2, 179--192.   

\bibitem{tosatti} Tosatti, Valentino: Limits of Calabi-Yau metrics 
when the Kahler class degenerates. ArXiv: 0710.4579. 

\bibitem{tsu} Tsuji, Hajime: Existence and degeneration 
of Kaehler-Einstein metrics on minimal algebraic varieties 
of general type. Math. Ann. 281(1988), no. 1, 123--133.

\bibitem{tsu2} Tsuji, Hajime: Degenerate 
Monge-Amp\`ere equation in algebraic geometry. 
Miniconference on Analysis and Applications 
(Brisbane, 1993), 209--224, Proc. Centre Math. 
Appl. Austral. Nat. Univ., 33, Austral. Nat. Univ., 
Canberra, 1994.

\bibitem{yau} Yau, Shing Tung: On the Ricci curvature of 
a compact Kaehler manifold and the complex Monge-Amp\`ere 
equation. I. Comm. Pure Appl. Math. 31(1978), no. 3, 339--411.

\bibitem{zh} Zhang, Zhou: On Degenerate Monge-Amp\`ere 
Equations over Closed K\"ahler Manifolds. Int. Math. Res. Not. 
2006, Art. ID 63640, 18 pp. 

\bibitem{thesis} Zhang, Zhou: Degenerate Monge-Amp\`ere 
Equations over Projective Manifolds. PHD Thesis at MIT, 2006. 

\bibitem{scalar-curvature} Zhang, Zhou: Scalar curvature bound 
for K\"ahler-Ricci flows over minimal manifolds of general type. 
ArXiv: 0801.3248. 

\end{thebibliography}
\end{document}